\documentclass[12pt]{amsart}
\usepackage{amsmath}
\usepackage{amsxtra}
\usepackage{amscd}
\usepackage{amsthm}
\usepackage{amsfonts}
\usepackage{amssymb}
\usepackage{mathrsfs}
\usepackage[all]{xy}
\usepackage[pdftex]{color,graphicx}
\usepackage{graphicx}
\usepackage{color}
\usepackage{mathrsfs}

\newtheorem{thm}{Theorem}

\newtheorem{ques}[thm]{Question}

\def\Z{\Bbb Z}

\def\Q{\Bbb Q}
\def\C{\Bbb C}

\def\R{\Bbb R}

 \DeclareMathOperator\SL{SL}

\makeatletter

\newcommand{\Rmnum}[1]{\expandafter\@slowromancap\romannumeral #1@}
\makeatother

\begin{document}

\title{Generating pairs for finite index subgroups of $SL(n,\Z)$}
\author{Chen Meiri}
\maketitle

\begin{abstract} Let $n \ge 3 $. Lubotzky \cite{Lu86} asked if every finite index subgroup $\Gamma$ of $\SL(n,\Z)$ contains a finite index subgroup which is generated by two elements. In this note we  show that for ``almost every'' element $g \in \SL(n,\Z)$ there exists $h \in \SL(n,\Z)$ such that $[\SL(n,\Z):\langle g,h^m\rangle] < \infty$ for every $m \ge 1$.  In particular, Lubotzky's question has a positive answer. 
\end{abstract}

\section{introduction}
Lubotzky \cite{Lu86} proved that if $\Gamma$ is a linear group which does not have a finite index meta-abelian subgroup then $\Gamma$
has a finite index subgroup  $\Lambda$ such that the minimal number of generators of every finite index normal subgroup of $\Lambda$ is at least three.  He then asked whether the restriction to normal subgroups is necessary. More specifically, he was interested in the case where $\Gamma=\SL(n,\Z)$. 

\begin{ques}[Lubotzky, \cite{Lu86}]\label{ques} Let $n \ge 3$. Does every finite index subgroup of $\SL(n,\Z)$ contain a finite index subgroup which is generated by two elements?\footnote{The answer is negative when $n=2$ since $\SL(2,\Z)$ contains a finite index  non-abelian free subgroup.}
\end{ques}

Previous works indicated that the answer for Question 1 should be positive:
\begin{enumerate}
\item Lubotzky--Mann,  \cite{LM87}: Let $p$ be a prime number and denote the $p$-adic integers by $\Z_p$. Every finite index subgroup of $\SL(n,\Z_p)$ contains a finite index subgroup which is generated by two elements.
\item Sharma--Venkataramana,  \cite{SR05}: Let  $\Gamma$ be a subgroup of finite index in $\mathrm{G}(\Z)$, where $\mathrm{G}$ is a connected semi-simple algebraic group over $\Q$ and of $\R$-rank $\ge 2$. If $\mathrm{G}$ has no connected normal subgroup defined over $\Q$ and $\mathrm{G}(\R)/\Gamma$ is not compact, then $\Gamma$ contains a subgroup of finite index generated by at most three elements.
\item Long--Reid,  \cite{LR11}: There exists a decreasing sequence $(\Lambda_i)_{i \ge 1}$
of finite index subgroups of $\SL(3,\Z)$ such that $\cap_{i\ge 1}\Lambda_i=1$ and $\Lambda_i$ is generated by two elements for every $i \ge 1$.
\end{enumerate}

Our main theorem is:

 \begin{thm}\label{main} Let $n\ge 3$ and $g \in \SL(n,\Z)$. Let $f_1,\ldots,f_r$ be the characteristic  polynomials of the blocks of the canonical rational form of $g$.  If the $f_i$'s are pairwise coprime then there exists a unipotent element  $u \in \SL(n,\Z)$ such that $[\SL(n,\Z):\langle g,u^m\rangle]<\infty$ for every $m \ge 1$.
 \end{thm} 
 
If the discriminant  of $g \in \SL(n,\Z)$ does not equal to zero then $g$ satisfies the assumption of Theorem \ref{main}. Every finite index subgroup of $\SL(n,\Z)$ is Zariski-dense in $\SL(n,\C)$ and must contain an element which does not belong to the subvariety  $\{g \in \SL(n,\C) \mid \mathrm{disc}(g)=0\}$. Thus, Theorem 2 provides a positive answer to Question 1.  
 
 The main ingredient in the proof of Theorem 2 is the congruence subgroup property or more presicely the following theorem of Bass, Lazard and Serre \cite{BLS64}: For every $n \ge 3$ and every $m \ge 2$ the subgroup $\langle e_{i,j}(m) \mid 1\le i \ne j \le n \rangle$ is equal to the $m$-congruence subgroup $\Gamma(m)$ of $\SL(n,\Z)$ where $e_{i,j}(m)$ is the matrix with 1 on the diagonal, $m$ on the $(i,j)$-entry and zero elsewhere while $\Gamma(m)$ consists of the matrices which are equal to the identity modulo $m$. Some of the  ideas used in the proof are similar to the ones in Venkataramana's paper \cite{Ve87}.
 
Our interest in Question 1 comes from the study of thin groups which are infinite index Zariski-dense subgroups of $\SL(n,\Z)$. 
This study is a very active field of research and yielded some remarkable theorems (see \cite{Sa14} and the reference therein). Yet, almost nothing is known about the algebraic structure of thin groups. The known examples of thin groups consist of  either free groups (or more generally certain free products) or surface groups and triangular groups \cite{LRT11}. 
It seems hard to construct other types of thin groups. 
One can explicitly calculate a subvariety of $\mathrm{V}$ of $\SL(n,\C)\times \SL(n,\C)$ such that if $(g,h)\in \SL(n,\C)\times \SL(n,\C) \setminus V$ then $\langle g,h \rangle $ is Zariski-dense. 
The main challenge is to understand when $[\SL(n,\Z):\langle g,h \rangle]<\infty$.  Theorem \ref{main} is a small step in this direction. 
\\ \\ 
{\bf Acknowledgement.} We are thankful to Nir Avni and Tsachik Gelander for several fruitful talks regarding similar questions. 

\section{proof of Theorem \ref{main}}

Assume that  $g \in \SL(n,\Z)$ satisfies the assumptions of Theorem \ref{main}.
There exists $v \in \Q^n$ such that $v,gv,\ldots, g^{n-1}v$ is a $\Q$-base of $\Q^n$ and there exists a $\Z$-base $w_1,\ldots w_n$ of 
$\Z^n$ such that for every $1 \le k \le n$, $$\mathrm{span}_\Q\{w_i \mid 1 \le i \le k\}=\mathrm{span}_\Q\{g^{i-1}v \mid 1 \le i \le k\}.$$ Thus, $g$ is conjugate in $\SL(n,\Z)$ to a matrix which belongs to
the Bruhat cell $\mathcal{B}_\sigma$ corresponding to the permutation $\sigma:=(12\cdots n)$.  Thus, we can assume that $g \in \mathcal{B}_\sigma$ and then it is enough to show that $[\SL(n,\Z):\langle g,e_{1,n}(m)\rangle]<\infty $ for every $m \ge 1$. Fix $m \ge 1$ and denote $\Lambda:=\langle g,e_{1,n}(m)\rangle$. Let $p_\sigma$ be the matrix with 1 on the $(\sigma(i),i)$-entry for every $1 \le i \le n-1$,
$(-1)^{n+1}$ on the $(\sigma(n),n)$-entry and zero elsewhere.

For every $1 \le i \ne j \le n$ let $E_{i,j}$ be the subgroup of $\SL(n,\C)$ consisting of the matrices with 1 on the diagonal and zero on the $(r,s)$-entry whenever $1 \le r \ne s \le n$ and $(r,s) \ne (i,j)$.  Moreover, for $1 \le i < j \le n$ define $$U_{i,j}:=\langle E_{r,s} \mid j<s \text{ or } j=s \text{ and } r \le  i \rangle$$ and 
$$V_{i,j}:=\left\{ \begin{array}{cc}
U_{i-1,j} & \text{if } i \ne 1 \\
U_{j,j+1} & \text{if } i = 1 \text{ and } j \ne n \\
\{\textrm{Id}\} & \text{if } (i,j)=(1,n)
\end{array}\right..$$

Since $g \in \mathcal{B}_\sigma$ there exist unique rational matrices $g \in B$ and $u \in U_{n-1,n}$ such that $g=bp_\sigma u$ where $B$ is the upper triangular subgroup of $\SL(n,\C)$. Let $2 \le i \le n-1$, since $U_{n-1,n}$ is abelian and $bE_{1,n}b^{-1}=E_{1,n}$
$$[gE_{i-1,n}g^{-1},E_{1,n}]=b[E_{i,1},E_{1,n}]b^{-1}=bE_{i,n}b^{-1}.$$ 
Thus, $\Lambda$ has a non-trivial intersection with $U_{i,n}\setminus V_{i,n}$. 
An induction argument implies that $\Lambda \cap U_{n-1,n}$ is Zariski-dense in $U_{n-1,n}$. 

For every $1 \le i < j \le n$ the group $U_{i,j}$ is normalized by $b$. In particular $b^{-1}U_{i,j}b$ contains 
$e_{i,j}(1)$.  Thus, 
if $2 \le i < j \le n$ then 
$p_\sigma^{-1}b^{-1}U_{i,j}bp_\sigma$ contains $e_{i-1,j-1}( 1)$ so $g^{-1}U_{i,j}g$ has a non-empty intersection with $U_{i-1,j-1}\setminus V_{i-1,j-1}$. An induction argument implies that 
$\Lambda \cap U$ is Zarisk-dense in $U=U_{1,2}$. Thus, $\Lambda$ contains $U_{1,2}(k\Z)$ for some $k \ge 1$. 

Since $U_{n-1,n}$ is abelian, $ge_{n-1,n}(k)g^{-1}=be_{n,1}(k)b^{-1}$. Denote $\Lambda^*:=b^{-1}\Lambda b$, then $e_{n,1}(k)$ belongs to $\Lambda^*$. Moreover, $U(d\Z) \subseteq b^{-1}
U_{n-1,n}(k\Z)b \subseteq\Lambda^*$ for some $d \ge 2$. For every $1 <i \le n-1$, $[e_{i,n}(d),e_{n,1}(k)]=e_{i,1}(kd)\in \Gamma$. Finally, if $1 \ne i \ne j \ne 1$ then $[e_{i,1}(kd),e_{1,j}(d)]=e_{i,j}(kd^2)\in\Lambda^*$. Thus, $\Lambda^*$ contains $\{e_{i,j}(kd^2) \mid 1 \le i \ne j \le n\}$. The proof of Bass, Lazard and Serre for the congruence subgroup property  \cite{BLS64} implies that $\Lambda^*$ contains the $kd^2$-congruence subgroup $\Gamma(kd^2)$ of $\SL(n,\Z)$.
Thus, $\Lambda^*\cap \SL(n,\Z)$ has a finite index in $\SL(n,\Z)$ so $\Lambda=b\Lambda^*b^{-1}$ also has a finite index in $\SL(n,\Z)$.

\end{document}